\documentclass[english]{article}
\usepackage[T1]{fontenc}
\usepackage[latin9]{inputenc}
\usepackage{amsmath}
\usepackage{graphicx}

\usepackage{babel}

\begin{document}

\title{On the Three Colorability of Planar Graphs}

\date{I. Cahit}
\maketitle
\begin{abstract}
The chromatic number of an planar graph is not greater than four and
this is known by the famous four color theorem and is equal to two
when the planar graph is bipartite. When the planar graph is even-triangulated
or all cycles are greater than three we know by the Heawood and the
Grotszch theorems that the chromatic number is three. There are many
conjectures and partial results on three colorability of planar graphs
when the graph has specific cycles lengths or cycles with three edges
(triangles) have special distance distributions. In this paper we
have given a new three colorability criteria for planar graphs that
can be considered as an generalization of the Heawood and the Grotszch
theorems with respect to the triangulation and cycles of length greater
than $\ge4$. We have shown that an triangulated planar graph with
$k$ disjoint holes is 3-colorable if and only if every hole satisfies
the parity symmetric property, where a hole is a cycle (face boundary)
of length greater than $3$. 
\end{abstract}

\section{Introduction}

Four coloring of planar graphs is a famous theorem (4CT) and it has
been proved twice by the same method by the assistance of an computer
and correctness of the proof has been verified by another computer
program {[}1],{[}2],{[}3]. The author has given an non-computer proof
of the four color theorem by using spiral chains and spiral chain
coloring in the maximal planar graphs {[}2],{[}3],{[}17]. However,
for an given planar graph the question {}``When three color suffice?''
has not been completely solved. When the planar graph is even-triangulated
or all cycles are greater than three we know by the Heawood and the
Grotszch theorems that the chromatic number is three {[}1]. In the
literature, there are several proofs of Grotszch theorems {[}11],{[}15],{[}16],{[}18]
and the simplest and efficient algorithmic proof is appear to be given
by the author {[}17]. It has been also considered on the other surfaces
{[}12],{[}13],{[}15]. Let $C=\{White,Gray,Black\}=\{W,G,B\}$ be the
set of three colors and when we needed we use $Red=\{R\}$ as the
fourth color.

In this paper we have given a new three colorability criteria for
planar graphs that can be considered as an generalization of the Heawood
and the Grotszch theorems with respect to the triangulation and cycles
of length greater than $\ge4$. First we have defined the triangulated
ring and gave necessary and sufficient condition for three colorability
ring. Next we have given an generaliztion of triangulated rings for
an triangulated planar graph with $k$ disjoint holes which is 3-colorable
if and only if every hole satisfies the parity symmetric property,
where a hole is a cycle (face boundary) of length greater than $3$.
\newpage{}

\begin{figure}
\includegraphics[angle=-90,origin=c,scale=0.4]{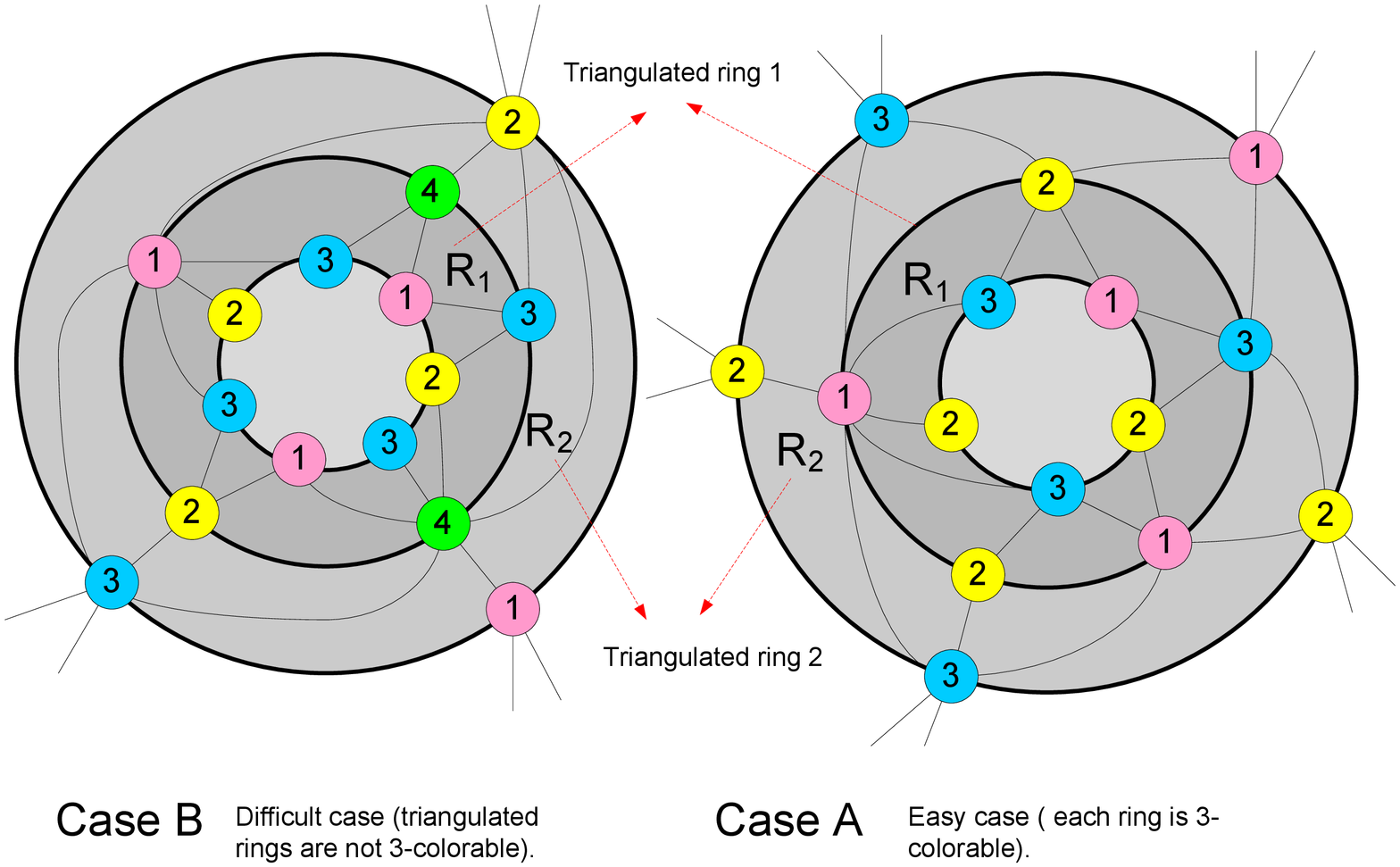}

\caption{Illustration of the triangulated rings for 3 and 4 colorings.}

\end{figure}
\begin{figure}
\includegraphics[scale=0.5]{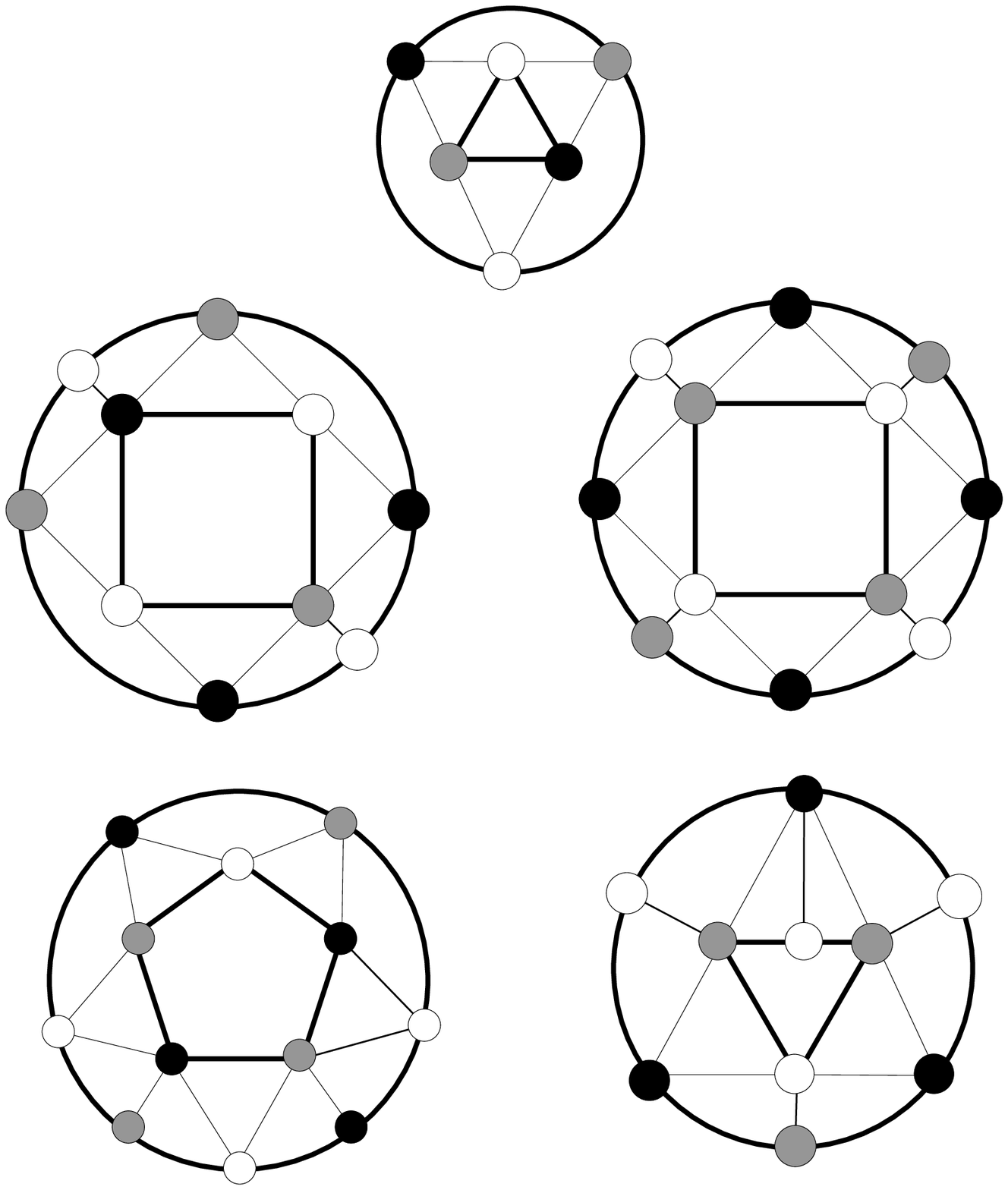}

\caption{All possible three colorings of the triangulated rings with inner-cycles
of length 3,4 and 5.}

\end{figure}

\begin{figure}
\includegraphics[scale=0.7]{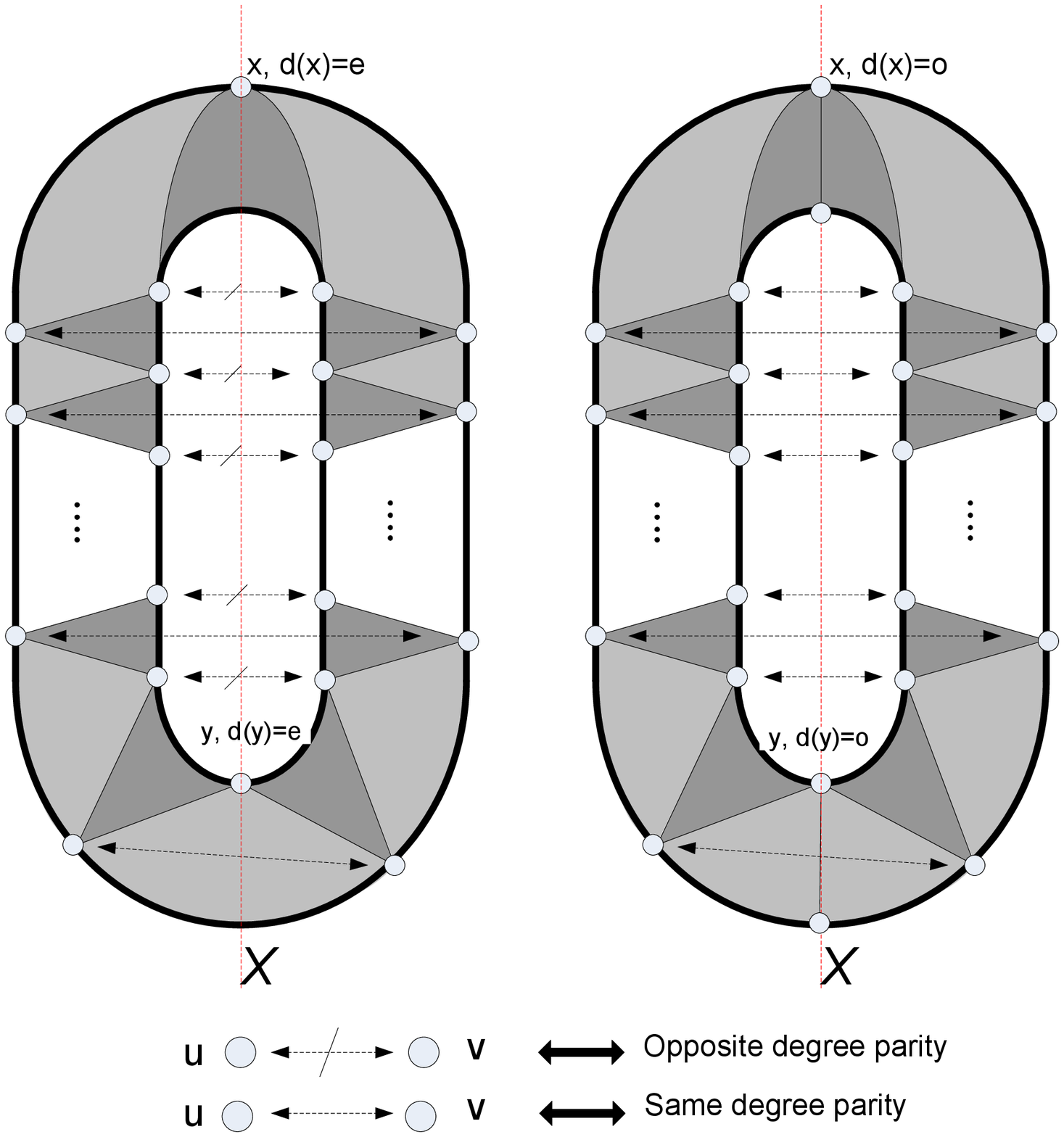}

\caption{Three colorable triangulated rings that $cps(G_{i})\in T,i=1,2$,
where $X$ is the symmetry axis of the graph.}

\end{figure}

\begin{figure}
\includegraphics[scale=0.6]{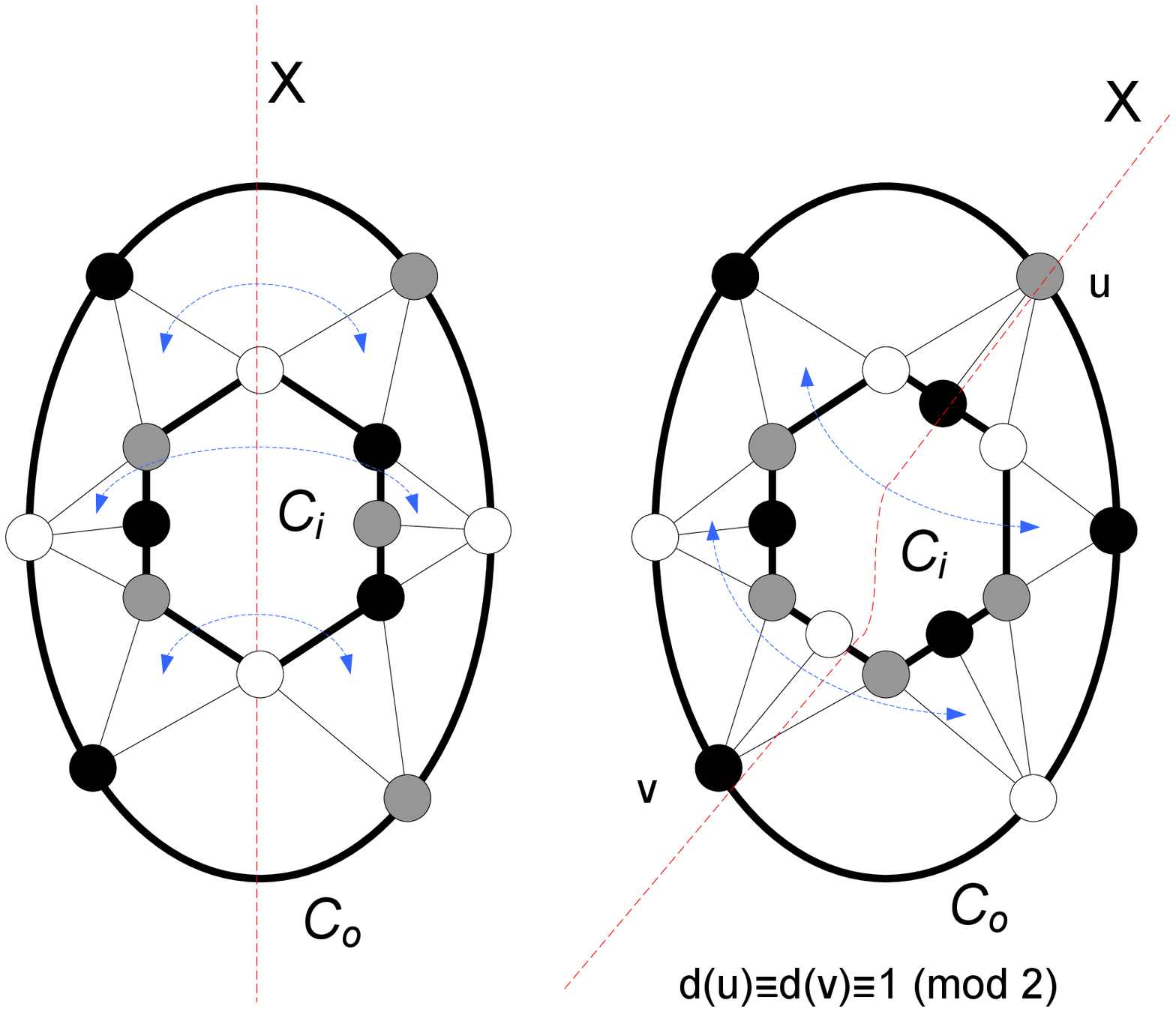}

\caption{Three colorings of the triangulated rings with $|C_{i}|=8,|C_{o}|=6$
and $|C_{i}|=10,|C_{o}|=6$}

\end{figure}

\section*{2 Triangulated ring}

In Fig. 1 we have shown two double triangulated rings together with
their colorings. The double triangulated ring shown on the right (Case
A) is $3$-colorable while double triangulated ring shown on the left
(Case B) is $4$-colorable.

We will make a few formal definitions first. A triangulated ring is
a $2$-connected planar graph $G$ with minimum degree $\geq3$ with
two faces $F_{i}$ and $F_{o}$ whose facial walks are the (induced)
cycles $C_{i}$ and $C_{o}$, respectively, such that:

(1) $V(C_{i})$ and $V(C_{o})$ partition $V(G)$. (That is $V(C_{i}){\cup}V(C_{o})=V(G)$
and  $V(C_{i}){\cap}V(C_{o})={\emptyset}$), where indices $i$ and
$o$ are being used to denote the inner and outer cycles (faces) of
the graph.

(2) Every face other than $F_{i}$ and $F_{o}$ is a triangle.

Every vertex $v$ on $C_{i}$ or $C_{o}$ has a fan, namely the triangles
incident with $v$ (other than possibly $C_{i}$ and $C_{o}$). These
fans, along with $F_{i}$ and $F_{o}$, partition the faces of $G$.
Let the fan graph of $G$ be the graph $F(G)$ whose vertices $f_{i}$
represent the fans $F_{i}$ of $G$, where an edge is in $F(G)$ if
the fans share a common edge. The vertices $f_{i}$ can (and will
be) identified with the of the fan (the vertex of the fan adjacent
to all others). The graph $F(G)$ is also the subgraph of $G$ induced
by $\{f_{0},f_{1}, ... , f_{2k-1}\}$. The fan graph of $G$ is a
even cycle $f_{0}, f_{1}, ... , f_{2k-1}$. The \emph{cyclic parity
sequence} of $(cps (G))$ is the cyclic sequence $p_{0}p_{1}p_{2} ... p_{2k-1}$
where is the parity (even/odd) of the number of triangles in $F_{i}$.
Note that every triangulated ring has at most one $3$-coloring (up
to permutation of colors), because a triangle has exactly one $3$-coloring
(again, up to permutation of colors), and every edge not in $C_{o}\cup C_{i}$
is in a triangle. The next lemma takes this result one step further.
For that, we must introduce some terminology related to the $cps$
of a graph. If $p_{0}p_{1} ... p_{2k-1}$ is a $cps$ with $k\ge 2$
and $p_{j} = e$, then define the $e$-collapse of the $cps$ at $j$
to be the $cps$ as $p_{0}p_{1} ... p_{j-2}(p_{j-1}+p_{j+1})p_{j+2} ... p_{2k}$,
where addition is of parities (modulo $2$). That is, the $e$-collapse
of $ooeeoeoo$ at $5$ is $ooee(o + o)o = ooeeeo$. Finally, let $T$
be the set of all $cps$'s which can be transformed into $ee$ or
$(o)^{6m}$ (i.e., $6m$$o$'s, for some integer $m$) by a finite
number of $e$-collapses.

Let us give a simple lemma.

\emph{Lemma 1}. \emph{If} $G(C_i{\cup} C_o)$ \emph{is 3-colorable
triangulated ring with all even fans then} $|C_i|\equiv 0 (mod 3)$\emph{or
with all odd fans then $|C_i|=|C_o|\equiv 0 (mod 2),|C_i|\not=4$.}

\emph{Lemma 2. If }$G$ \emph{is a triangulated ring, then} $G$ \emph{is
3-colorable iff} $cps (G)\in T$ . 

First, we need some lemmas. The first is easily proven using induction:

\emph{Lemma 2.1. Suppose }$F$\emph{ is a fan graph, where} $v$ \emph{is
the vertex of} $F$ \emph{adjacent to all others, and }$u$ \emph{and}
$w$ \emph{are the two vertices of $F\\v$ with degree $1$. Then:}

\emph{(a) If $F$ is even and $u,v$and $w$ are colored so that $u$
and $w$ receive the same color, then this partial $3$-coloring extends
to a proper $3$-coloring of $F$;}

\emph{(b) In a proper $3$-coloring of an even fan $F$, $u$ and
$w$ receive the same color;}

\emph{(c) If $F$ is odd and $u$, $v$, and w are colored so that
$u, v$ and $w$ all receive different colors, then this partial $3$-coloring
extends to a proper $3$-coloring of $F$; and}

\emph{(d) In a proper $3$-coloring of an odd fan $F, u, v$, and
$w$ all receive different colors.}

If $G$ is a triangulated ring with at least four fans, and $F_{j}$
an even fan of $G$, we will call the {[}even] fan collapse of $G$
at $F_{j}$ the graph $H$ obtained from $G$ by deleting the vertices
of $F_{j}$ other than $f_{j-1},f_{j}$ , and $f_{j+1}$, and then
identifying the vertices $f_{j-1}$ and $f_{j+1}$.

\emph{Lemma 2.2. Let} $G$ \emph{be a triangulated ring, and} $H$
\emph{the graph resulting from a fan collapse of} $G$ \emph{at} $F_{j}$
. \emph{Then} $G$ \emph{is} $3$-\emph{colorable if and only if }$H$\emph{
is. Furthermore,} $cps (H)$ \emph{can be obtained from} $cps (G)$
\emph{by an} $e$-\emph{collapse at}$j$.

\emph{Proof:} If $c$ is a proper $3$-coloring of $G$, then since
$F$ is a fan, $c(u) = c(w)$, by Lemma 2.1(b). This means that when
$H$ is created, some vertices will be deleted (which makes this smaller
graph $3$-colorable), and two vertices with the same color will be
identified. This means $H$ is $3$-colorable. Now suppose $c$ is
a proper $3$-coloring of $H$. To go from $H$ to $G$, a vertex
$v_{0}$ will have to be split into $u$ and $w$. We will let $c(u)= c(w) = c(v_{0})$,
and we still have a $3$-coloring. To get $G$, we need to add some
vertices, so we need to make sure that we can extend the $3$-coloring
to these vertices. Since the fan which will be created is even, Lemma
2.1(a) implies that this new fan is $3$-colorable. The union of these
two proper $3$-colorings gives a proper $3$-coloring of $G$.

To prove the second result: The graph $H$ is a triangulated ring
with two fewer fans than $G$. All but one of the fans of $H$ come
from fans of $G$, with the other fan of $H$ being the combination
of two fans $F_{j-1}$ and $F_{j+1}$ of G; the parity of this new
fan is the sum of the parities of $F_{j-1}$ and $F_{j+1}$. QED.

And the following result makes life easier:

Lemma 2.3. If $G_{1}$ and $G_{2}$ are triangulated rings with the
same $cps$, then $G_{1}$ is $3$-colorable iff $G_{2}$ is.

Proof: If $G_{1}$ and $G_{2}$ are as stated above, then $G_{1}$
can be transformed into $G_{2}$ by repeatedly adding two vertices
to a fan, or by removing two vertices from a fan. If this is done
to fan Fi, it doesn't affect the coloring of $f_{i-1}, f_{i}$, or
$f_{i+1}$. QED.

\emph{Proof of Lemma 2}: Let $F(G)$ be the fan graph of $G$. Note
that fan $F_{i}$ contains the vertices $f_{i_1}, f_{i}$, and $f_{i+1}$
(where all indices are modulo $2k$), which correspond to the vertices
$u, v$, and $w$ in the statement of Lemma 2.1.

We will begin by proving that if $cps (G)\in T$ , then $G$ is $3$-colorable.
The proof will be by induction on the number of $e$-collapses. If
there are no $e$-collapses, then we only have to settle the cases
$ee$ and $(o)^{6m}$, since $eo,oo\notin T$ . If $cps(G) = ee$,
then $F(G)$ contains two vertices $v_{1}$ and $v_{2}$, and (technically)
$2$ parallel edges. (Note that $G$ itself can itself be a simple
graph, because the $e$ in the $cps$ means there's an even number
of triangles in each fan, not zero.) If we color $f_{0}$ with $1$
and $f_{1}$ with $2$, then note that for all $i, p_{i} = e$, and
$f_{i-1}$ and $f_{i+1}$ are (trivially) colored with the same color.
Lemma 2.1(a) implies that $G$ is $3$-colorable. Similarly, if $cps (G)$
consists of $6m o$'s, we color the fan graph as follows. (Note that
in this case, $2k = 6m$ is a multiple of $6$, so it's also a multiple
of $3$, and the coloring is well-defined:

$c(f_{i})=\begin{cases}
1 & \text{if }i\equiv0(mod3)\\
2 & \text{if }i\equiv1(mod3)\\
3 & \text{if }i\equiv2(mod3)\end{cases}$

Now, since every $p_{i}$ is odd, all we need to do is to verify that
$f_{i-1}, f_{i}$, and $f_{i+1}$ all receive different colors, for
all $i$. This follows immediately from the definition, and Lemma
2.1(c) then implies that $G$ is $3$-colorable. Now suppose the result
is true for $N -1$ $e$-collapses, with $N \ge 1$. Let G be a triangulated
ring whose $cps (G)$ can be transformed into $ee$ or $(o)^{6m}$
by $N$ $e$-collapses. Consider the first $e$-collapse, which we
will assume occurs at $j$ and results in the $cps$ sequence $S$.
If $H$ is the fan collapse of $G$ at $F_{j}$ , then $cps (H) = S$,
by Lemma 2.2.

But since $S$ can be transformed into $ee$ or $(o)^{6m}$ with $N-1$
$e$-collapses, $H$ is $3$-colorable by the induction hypothesis.
Lemma 2.2 states that $H$ is $3$-colorable iff $G$ is, so $G$
is $3$-colorable.

Now we have to show that if a triangulated ring $G$ is $3$-colorable,
then $cps (G)\in T$ . Assume that $G$ is $3$-colorable. First,
we will settle the case where $cps (G)$ has no $e$'s in it. In this
case, every $p_{i}$ is $o$, so the$3$-coloring must satisfy $c(f_{i+3}) = c(f_{i})$
for all $i$, where indices are taken modulo $2k$; this follows because
$c(f_{i}), c(f_{i+1})$, and $c(f_{i+2})$ must all be distinct, since
$p_{i+1} = o$, by Lemma 4.1; and $c(f_{i+1}), c(f_{i+2})$, and $c(f_{i+3})$
must all be distinct, since $p_{i+2} = o$. This forces $c(f_{i}) = c(f_{i+3})$.
However, if the number of $o$'s is not a multiple of $6$, then $2k$
(the length of the $cps$) is not a multiple of $3$. That means that
the condition $c(f_{i+3}) = c(f_{i})$ forces all vertices $f_{i}$
to be colored the same, which is not allowed by Lemma 2.1(c). Hence
there is no proper $3$-coloring of $G$, contrary to assumption.
Thus the number of $o$'s is a multiple of $6$, so $cps (G) = (o)^{6m}\in T$
, as claimed. Now suppose that $p_{j} = e$. If $k = 1$, then $cps (G) = oe, eo$,
or$ee$. It is easily seen that if $cps(G) = oe$ (or $eo$) then
$G$ is not $3$-colorable, so $cps (G) = ee \in T$ .

So now we may assume that $k \ge 2$. Define a sequence of graphs
$G_{i}$ in the following way: Let $G_{0}$ be $G$, and for all $i \ge 0$,
if $p_{j} = e$ in $cps (G_{i})$ and $cps (G_{i})$ has length at
least two, then let $G_{i+1}$ be the fan collapse of $G_{i}$ at
$F_{j}$ . Suppose we cannot continue from $G_{N}$. Note that, for
all applicable $i$, $cps (G_{i+1})$ can be obtained from $cps (G_{i})$
by $e$-collapse. Then $cps (G_{N})$ either has no $e$'s, or has
length two. Furthermore, $G_{N}$ is $3$-colorable, by repeated application
of Lemma 2.2. But we have seen $G_{N}$ can only be $3$-colorable
if $cps (G_{N}) = ee$ or $cps (G_{N}) = (o)^{6m}$. In either case,
we can obtain $ee$ or $(o)^{6m}$ from $cps (G)$ by repeated $e$-collapsings.
Thus $cps (G)\in T$, which proves the lemma.

The following lemma is given without proof and useful for $3$-colorable
triangulated rings.

\emph{Lemma 3. Triangulated ring $G$ is $3$-colorable if $cps(G)$
is symmetric and $|C_{o}|\equiv0(mod3)$ or $|C_{o}\cup C_{i}|\equiv0(mod3)$.}

For illustrations see the $3$-colorings of the triangulated rings
given in Fig. 2 and 4. In Fig. 4 the symmetry axis is denoted by $X$.

In Fig. 3 a more general three colorable triangulated rings have been
shown. The $csp(G)$'s respectively are:

$cps(G_{1})=\{oab^{,}...ab^{,}ab^{,}obaba...ba\}$ (Fig.2 (left)) 

$cps(G_{2})=\{eab...babaebaba...ba\}$ (Fig.2 (right)),

where if $b=e$ then $b^{,}=0$ and if $b^{,}=e$ then $b=o$, and
$o$ and $e$ denote odd and even parities.

\section*{3 Planar graphs with holes}

In this section we will extend the result obtained for the three-colorability
of the triangulated rings to triangulated planar graphs with vertex
disjoint $k$ cycles (holes) $h_{i}$ of lengths greater than three.
Let $G(h_{i}),$ $i=1,2,...,k$ be the set of triangulated rings arround
the holes $h_{i}$, where $V(h_{1})\cap V(h_{2})\cap...\cap V(h_{k})=\emptyset$.
Let $V(C_{o})$ be the outer-cycle of $G$. Clearly $V(C_{o})\cap\{V(h_{1})\cup V(h_{2})\cup...\cup V(h_{k})\}=\emptyset$.

\emph{Theorem 3. Let $G$ be an triangulated planar graph with $k$
disjoint holes $h_{i}$, $i=1,2,...,k$. Then $G$ is $3$-colorable
iff for every triangulated ring $G(h_{i})$ we have $csp(G(h_{i}))\in T$.}

\emph{Proof. }Necessity of the theorem can be easily seen by the cyclic
parity sequence of $(cps G(h_{i})$ which we have assumed that can
be transformed into $ee$ or $(o)^{^{\{6m\}}}$(see Lemma 2). 

Now define the graph $H(V,E)$ where the vertex set $V(H)=\{h_{1},h_{2},...,h_{k}\}$
are the holes of $G$ and $e_{i,j}=(h_{i}h_{j})\in E(H)$ if $E(G((h_{i}))\cap E(G(h_{j}))\neq\emptyset$.
Let $T_{G}$ be an spanning tree of $H(V,E).$ Then start coloring
of the vertices $G(h_{1})$ first and next select an vertex $h_{i},i\neq1$such
that $(h_{i}h_{1})\in T_{G}.$ Color the vertices of $G(h_{\imath})$.
Repeat this step for the other vertices of $T_{G}$. It is clear that
since $cspG(h_{i})\in T$ and $(G(h_{i})$is an triangulated ring
at the end $G$ would be colored properly with three colors. 

\emph{Corollary}. The planar triangulated graph $G(h_{i})$ with $k$
holes $h_{i}=1,2,...,k$ can be made $3$-colorable triangulated ring
with $C_{i}=h_{1}\cup h_{2}\cup...\cup h_{k}$.

Proof. Delete suitable edges of the spanning tree $T_{G}$ in merging
two adjacent holes $h_{i}$ and $h_{j}$. At the end the inner cycle
$C_{i}$ will be the union of the holes. 

Next define the semi-triangulated graph $G^{*}$as an triangulated
ring in which there exists at least one cycle $C_{r}$ of length greater
than three such that $C_{i}\cap C_{r}\neq\emptyset$ and $C_{o}\cap C_{r}\ne\emptyset$.
Note that non-triangulated ring has no triangle; hence is $3$-colorable
by Grotzsch theorem. The following simple theorem gives useful information
when $G^{*}$is $3$-colorable.

\emph{Theorem 4. Let $G$ be an triangulated ring. If $csp(G)\notin T$
and $G^{*}$be any semi-triangulated ring obtained by the addition
of an single cycle $C_{r}$ of length $k\geq4$ then $G^{*}$is $3$-colorable.
If $csp(G)\in T$ and $G^{*}$be any semi-triangulated ring obtained
by the addition of an single cycle $C_{r}$ such that $|E(C_{o})\cap E(C_{r})|\ge2$
then $G^{^{*}}$is $3$-colorable.}

\emph{Proof. }If $csp(G)\notin G$ then for any $3$-coloring $c$
there must be an vertex $v$ such that $c(v)=i=j,$ where $i,j\in\{1,2,3\},i\neq j$.
Note that such a vertex $v$ can be selected freely beforehand since
$G$ is an triangulated ring. Then the vertex $v$ is splited into
two $v^{'}$and $v^{''}$vertices in $G^{*}$such that $(v^{'}v^{''})\in C_{r}$.
Hence the $3$-coloring $c$ of $G^{*}$can be made proper by $c(v^{'})=i$
and $c(v^{''})=j$ and the other vertices of $C_{r}$can be colored
alternatingly by two colors $1,2$ or $3,2.$ Second part of the theorem
is similar to the first part but this time since $csp(G)\in G$ the
coloring $c$ of $G$ is an proper $3$-coloring and the vertex $v$
with $c(v)=i$ must be splited into three vertices i.e., \emph{$|E(C_{o})\cap E(C_{r})|\ge2$.
}Therefore coloring $c$ will be proper three coloring of $G^{*}=G\cup C_{r}$.

In fact the above theorem can be generalized to semi-triangulated
rings with $k$ cycles of length $\ge4$. That is think of $G^{*}$as
disjoint triangular ladders (all with triangles) separated by cycles
of length $\ge4$. Let us denote the semi-triangulated graph $G^{*}=G^{*}(C_{1}\cup C_{2}\cup...\cup C_{k})$.
This suggest that as we are inserting large cycles into $3$-colorable
planar graph three colorability maintained. 

\begin{figure}
\includegraphics[scale=0.5]{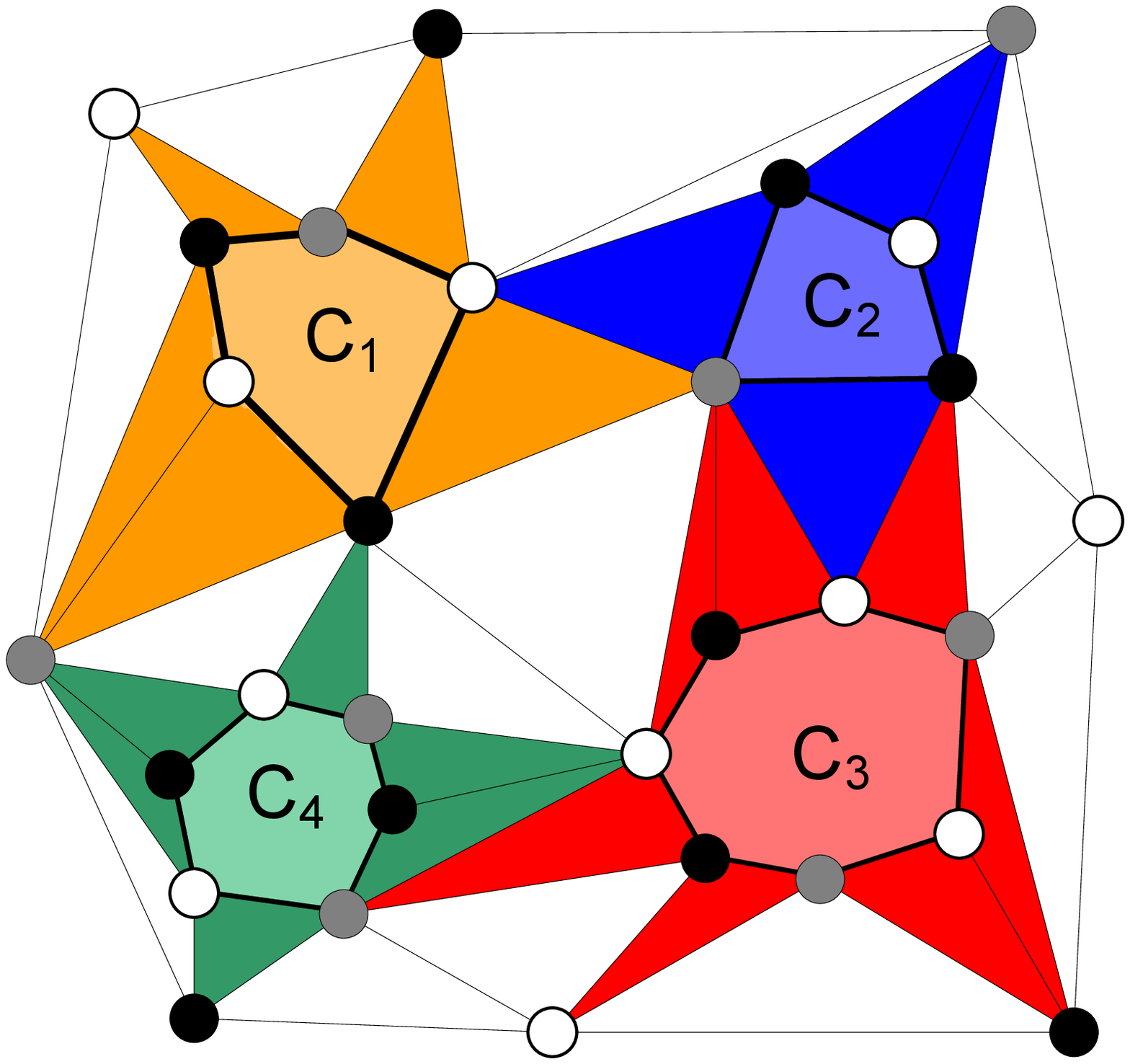}

\caption{Three coloring of an planar graph with $4$ holes.}

\end{figure}

\begin{figure}
\includegraphics[scale=0.6]{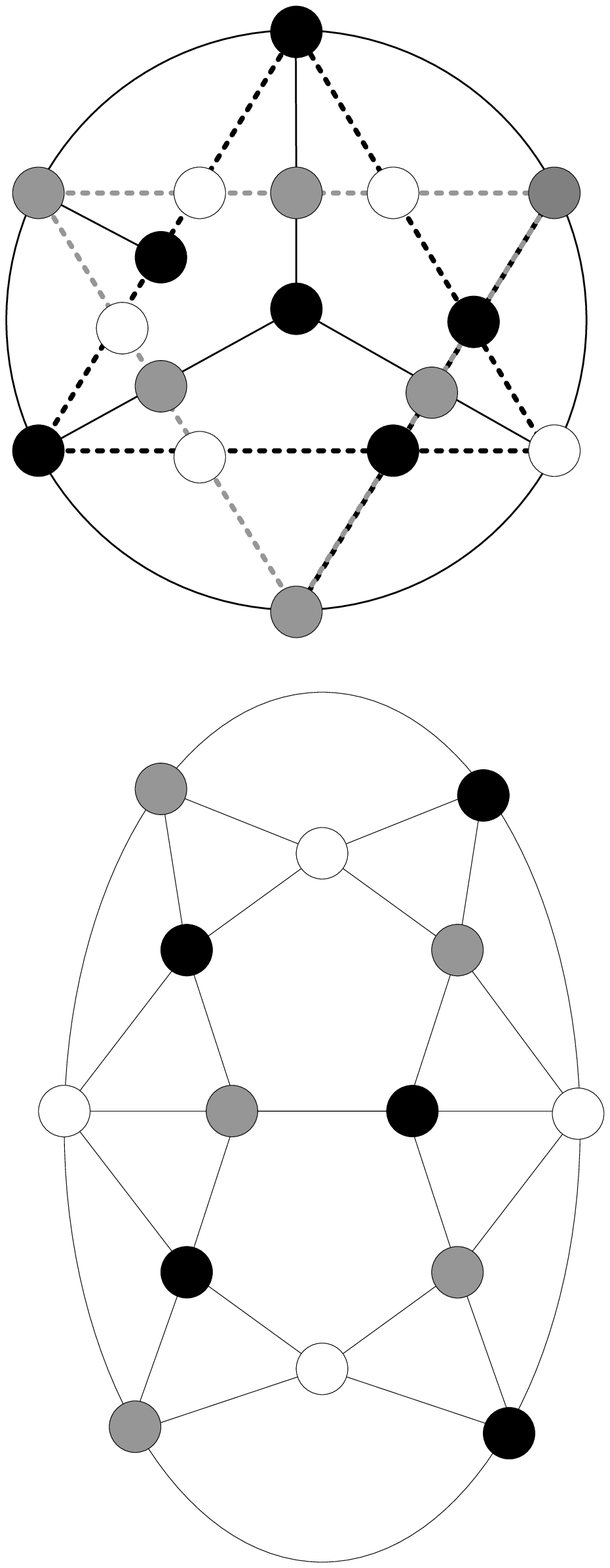}

\caption{Three colorable planar graphs with triangulated ring (\emph{upper}):
with three cycles of length 5 and 6, (\emph{lower}): with two cycles
of length 5.}

\end{figure}

\section*{4 Another wave of conjectures on $3$-colorability}

As early as 1959, Grötzsch proved that every planar graph without
$3$-cycles is $3$-colorable. This result was later improved by Aksionov
in 1974. He proved that every planar graph with at most three $3$-cycles
is $3$-colorable. In 1976, Steinberg conjectured the following :

\emph{Steinberg's Conjecture (1976) {[}5]} : Every planar graph without
$4$- and $5$-cycles is $3$-colorable. 

An algorithmic proof to Steinberg's conjecture has been proposed by
the author in 2006 {[}4]. We note that the statement of the Steinberg's
conjecture is not sharp since there are $3$-colorable planar graphs
with four and five cycles (see Fig. 7). 

In 1969, Havel posed the following problem:

\emph{Havel's Problem (1969) {[}20]} : Does there exist a constant
$C$ such that every planar graph with the minimal distance between
triangles at least $C$ is $3$-colorable?

Aksionov and Mel'nikov proved that if $C$ exists, then $C\ge4$,
and conjectured that $C=5$ {[}6]. These two problems remain widely
open. In 1991, Erdös suggested the following relaxation of Steinberg's
Conjecture: Determine the smallest value of k, if it exists, such
that every planar graph without any cycles of length $4$ to $k$
is 3-colorable. The best known bound for such a $k$ is $7$ {[}19].
Many other sufficient conditions of $3$-colorability considering
planar graphs without cycles of specific lengths were proposed.

At the crossroad of Havel's and Steinberg's problems (since the authors
consider planar graphs without cycles of specific length and without
close triangles), Borodin and Raspaud proved that every planar graph
without $3$-cycles at distance less than four and without $5$-cycles
is $3$-colorable (the distance was later decreased to three by Xu,
and to two by Borodin and Glebov). As well, they proposed the following
conjecture:

\emph{Strong Bordeaux Conjecture (2003)} {[}21]: Every planar graph
without $5$-cycle and without adjacent triangle is $3$ colorable.

By adjacent cycles, we mean those with an edge in common. This conjecture
implies Steinberg's Conjecture. Finally, Borodin et al. considered
the adjacency between cycles in planar graphs where all lengths of
cycles are authorized, which seems to be the closest from Havel's
problem ; they proved that every planar graphs without triangles adjacent
to cycles of length from $3$ to $9$ is $3$-colorable . Moreover
they proposed the following conjecture:

\emph{Novosibirsk 3-Color Conjecture (2006)} {[}22]: Every planar
graph without $3$-cycles adjacent to cycles of length $3$ or $5$
is $3$-colorable. This implies Strong Bordeaux Conjecture and Steinberg's
Conjecture.

We claim that the above two conjectures can be settled by the use
of spiral chains and this will be given in the updated version of
{[}4].

In Fig. 9 we have shown triangulated rings with $3$- and $4$-colorings.
Note that the Kempe chains $K(W,B),K(G,B)$ in Fig. 9(b) and $K(W,G)$
in Fig. 9(c) prevent the three colorability of the triangulated rings.

\section*{4 Concluding remarks}

In this paper we have given a new three coloring criteria, that is
three colorability of triangulated rings. We have generalized three
colorable rings to some extend to the triangulated planar graphs with
disjoint large cycles (cycles of length $\geq4$). Although we know
that planar $3$-colorability is \emph{NP}-complete {[}10] but the
result obtained here may rise some hopes to devise an efficient algorithm
for three colorable planar graphs. 

\smallskip{}
\smallskip{}

\begin{figure}
\includegraphics[scale=0.7]{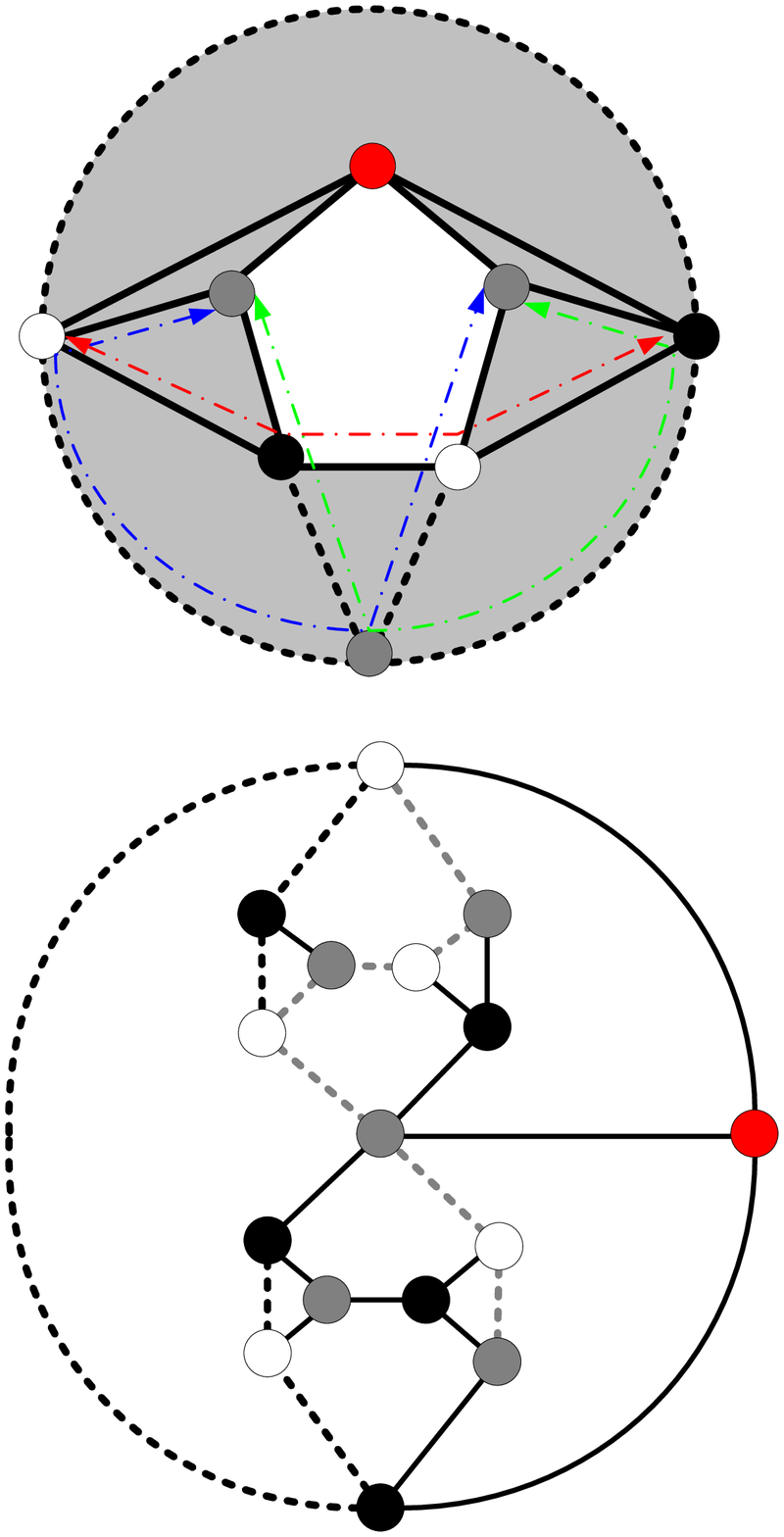}

\caption{(a) Four colorable triangulated ring with only an single $C_{5}$
(subgraph in bold lines is not $3$-colorable).(b) Four colorable
graph without $C4$.}

\end{figure}

\begin{figure}
\includegraphics[scale=0.6]{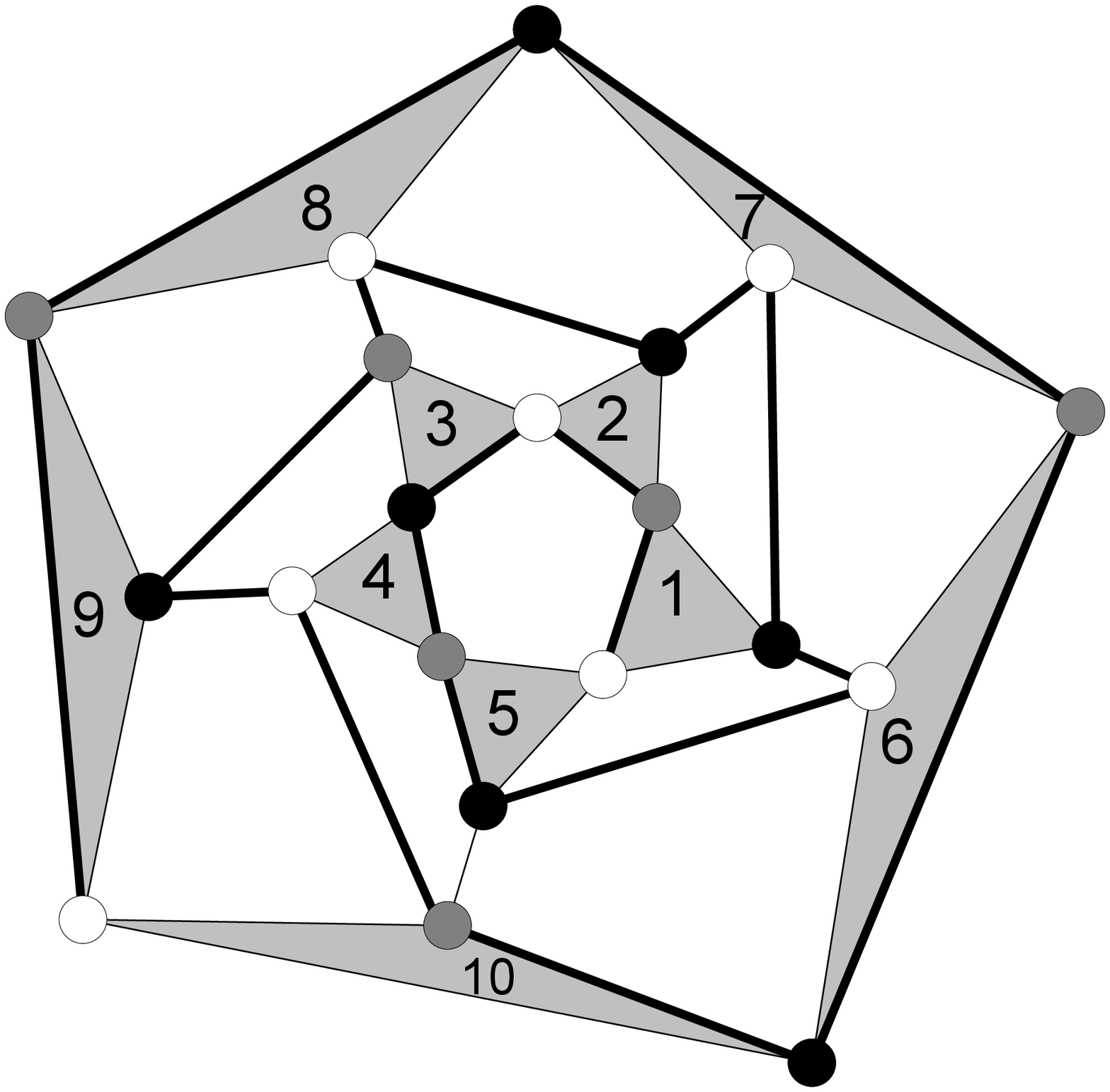}

\caption{Three colorable semi-triangulated rings.}

\end{figure}

\author{\textbf{Acknowledgments}. {\normalsize The author would like to thank
to Dr. Chris  Heckman for introducing {}``collapsing'' in the cps
sequences and the proof of Lemma 2.}}

\end{document}